\documentclass[12pt]{article}
\usepackage{amsmath}
\parskip=\medskipamount
\newtheorem{definition}{Definition}[section]

\newtheorem{theorem}[definition]{Theorem}
\newtheorem{proposition}[definition]{Proposition}

\newtheorem{remark}[definition]{Remark}

 1   
\font\ddpp=msbm10  scaled \magstep 1  
 1  
\def\bull{\ \ \ \vrule height 1.5ex width.8ex depth.3ex \medskip}

\def\QED{\hskip0.1em\hfill\null\ \null\nobreak\hfill
\kern3pt\lower1.8pt\vbox{\hrule\hbox   {\vrule\kern1pt\vbox{\kern1.7pt
\hbox{$\scriptstyle   QED$}\kern0.2pt}\kern1pt\vrule}\hrule}}
\def\R{\hbox{\ddpp R}}    

\def\ene{\hbox{\ddpp N}}    

\def\appendix{\par
 \setcounter{section}{0}
 \setcounter{subsection}{0}
 \def\thesection{Appendix \Alph{section}}}
\title{\bf Discrete variational integrators and optimal control theory}

\author{Manuel DE LE\'ON
\thanks{mdeleon@imaff.cfmac.csic.es}\and David MART\'IN DE DIEGO\thanks{d.martin@imaff.cfmac.csic.es}\and  Aitor SANTAMAR\'IA-MERINO\thanks{aitors@imaff.cfmac.csic.es}
\bigskip
\\
Laboratory of Dynamical Systems, Mechanics and Control\\
Instituto de Matem\'aticas y F{\'\i}sica Fundamental\\
Consejo Superior de Investigaciones Cient{\'\i}ficas\\
Serrano 123, 28006 Madrid, SPAIN }

\begin{document}
\maketitle
\abstract{A geometric derivation of  numerical integrators for optimal control problems is proposed.
It is based in the classical technique of generating functions adapted to the special features of  optimal control problems. }

\section{Introduction}

Optimal control has been one of the driving forces behind many of the applications of mathematics to 
engineering, robotics, economics... In fact, the Maximum Principle was discovered by L.S. Pontryagin in 1955 in an attempt to find a solution for a highly specific optimization problem related to the manoeuvres of an aircraft.
 One of its main features is the interplay among different research areas, specially control theory, classical mechanics and differential geometry. Historically,  Optimal Control Theory (OCT) took place during the 1950's and its geometrization   was started in the 1960's. This geometric analysis of OCT has been introduced using many fundamental tools of differential geometry: Lie groups, exterior differential systems, fiber bundles, riemannian and subriemannian geometry among others.

From other point of view, a geometric  methodology has been  recently shown to be very useful  
 for  simulating numerically the motion of  dynamical systems.
Following this research line, new numerical methods have been  developed, called geometric integrators;  usually, these integrators, in simulations, can  run for longer times with lower spurious  effects (for instance, bad energy behavior for conservative systems) than  the traditional (non-geometrical) ones.  
In particular, we are interested in extensions to OCT of Discrete variational integrators.
 These integrators have  precisely   their roots in the optimal control literature  in the 1960's  and 1970's (Jordan and Polack \cite{Jordan}, Cadzow \cite{Cadz}, Maeda \cite{Maed1,Maed2}) and in 1980's by Lee \cite{Lee1,Lee2}, Moser and  Veselov \cite{Mose}. Although this kind of symplectic integrators have been considered  for conservative systems 
\cite{Jaro1,Kane1,Mars6},   it  has been recently shown how  discrete variational mechanics can include forced or dissipative systems \cite{Kane3,Mars6}, holonomic constraints \cite{Mars6},  time-dependent systems \cite{LD1,Mars6}, frictional contact \cite{Pand} and nonholonomic constraints (see \cite{Cort,JS,LDA,LDA1}). Moreover, it has been also discussed reduction theory \cite{Bobe,Bobe1},   extension to field theories \cite{Jaro2,Mars2} and quantum mechanics \cite{Nort}.

In this paper, we shall continue this work by extending to the discrete variational techniques to Optimal Control Problems and relating our results with Discrete Optimal Control Theory. Mainly, we shall give a geometrical construction of symplectic integrators for OCT, proving as a direct consequence the symplecticity of some discrete optimal control problems. As a nice consequence,  an easy proof of the symplecticity of discrete Hamilton equations will be given.

Since most engineering systems are time-dependent, we shall include the time variable explicitely in our control models
and some geometric tools (mainly, cosymplectic geometry)
of time-dependent mechanics will be useful 

\section{Optimal control theory}\label{optimal}

It is well known that the dynamics of a large class of  engineering and  economic systems can be expressed as a set of differential equations
\begin{equation}\label{general}
\dot{q}^A=\Gamma^A(t, q(t),u(t)) \,,\; 1 \le A \le n \, ,
\end{equation}
where $t$ is the time,  $q^A$ denote the state variables and  $u^a$, $1\leq a\leq m$, the control inputs to the system that must be specified.
Given an initial condition of the state variables and given control inputs we completely know the trajectory of the state variables $q(t)$ (all the functions are assumed to be at least $C^2$).

 Given an initial condition, usually $q_0=q(t_0)$,  our aim  is to find a $C^2$-piecewise smooth curve $\gamma(t)=(q(t),u(t))$, satisfying the control equations 
(\ref{general}) and minimizing the functional 
\begin{equation}\label{general1}
{\mathcal J}(\gamma) = \int^{T}_{t_0} L(t, q(t),u(t))\,dt+S(T, q(T)) \, ,
\end{equation}
for some fixed and given final time $T\in \R^+$. The  integral  $\int^{T}_{t_0} L(t, q(t),u(t))\,dt$ depends on the time history (from $t_0$ to $T$) of the state variables and the control inputs, and  $S(\cdot, q(\cdot))$ is a cost function based on the final time and the final states of the system. 

In a global description, one assumes a fiber bundle structure $\pi: \R\times C 
\longrightarrow Q$, where $Q$ is the configuration manifold with local 
coordinates $(q^A)$ and $C$ is the bundle of controls, with  
coordinates $(q^A, u^a)$, $1\leq A\leq n$, $1\leq a\leq m$.

The time-dependent ordinary differential equations (\ref{general}) on $Q$ depending on 
the parameters $u$ can be seen as a vector field $\Gamma$ along the 
projection map $\pi$, that is, $\Gamma$ is a smooth map 
$\Gamma: \R\times C \longrightarrow TQ$ such that the diagram
\begin{figure}[h]
\centering
\setlength{\unitlength}{1cm}
\begin{picture}(4,2.5)(-0.7,0)
\put(0,2){\makebox(0,0)[r]{$\R\times C$}}
\put(4,2){\makebox(0,0)[l]{$TQ$}}
\put(2,0){\makebox(0,0)[c]{$Q$}}
\put(0.2,2){\vector(1,0){3.6}}
\put(2,2.3){\makebox(0,0)[r]{$\Gamma$}}
\put(0.2,1.8){\vector(1,-1){1.6}}
\put(3.8,1.8){\vector(-1,-1){1.6}}
\put(0.8,0.8){\makebox(0,0)[r]{$\pi$}}
\put(3.2,0.8){\makebox(0,0)[l]{$\tau_Q$}}
\end{picture}
\vspace{-0.5cm}
\centering
\end{figure}

\noindent is commutative. This vector field is locally written as 
$\Gamma=\displaystyle{\Gamma^A(t,q,u) \frac{\partial}{\partial q^A}}$.

A neccesary condition for the solutions of such problem are provided by Pontryaguin's maximum principle. 
If we construct the pseudo-Hamiltonian function:
\begin{equation}\label{hamiltonian}
H(t, q, p, u)=p_A \Gamma^A(t,q, u)-L(t, q, u)=p\Gamma(t, q,u)-L(t, q,u)
\end{equation}
where $p_A$, $1\leq A\leq n$,  are now considered as Lagrange's multipliers, then a curve $\gamma: [t_0, T]\rightarrow  C$, $\gamma(t)=(q(t), u(t))$ is an optimal trajectory  if there exist functions  $p_A(t)$, $1\leq A\leq n$,  such that they are  solutions of the pseudo-Hamilton equations:
\begin{equation}\label{eqH}
\left\{\begin{array}{l}
\dot{q}^A(t)=\displaystyle{\frac{\partial H}{\partial p_A}(t, q(t), p(t), u(t))}\\
\,\\
\dot{p}_A(t)=-\displaystyle{\frac{\partial H}{\partial q^A}(t, q(t), p(t), u(t))}\\
\end{array}\right.
\end{equation}
and we have
\begin{equation}\label{H}
H(t, q(t), p(t), u(t))=\underset{v}{\hbox{min}}\; H(t, q(t), p(t), v), \quad t\in [t_0, T]
\end{equation}
with transversality conditions
\[
q(0)=q_0\qquad \hbox{and}\qquad p_A(T)=-\frac{\partial S}{\partial q^A}(T, x(T))
\]

Condition (\ref{H}) is usually replaced by 
\begin{equation}\label{constraint}
\frac{\partial H}{\partial u^a}=0, \quad 1\leq a\leq m\; ,
\end{equation}
when we are looking for extremal trajectories.

It is well known that the Pontryaguin's necessary conditions for extremality have a geometric interpretation 
in terms of presymplectic (or precosymplectic) Hamiltonian systems. 
The total space of the system will be $\R\times (T^*Q \times_Q C)$, with induced coordinates $(t, q^A, p_A, u^a)$ .

Define the Pontryaguin's Hamiltonian function $H: \R \times (T^*Q \times_Q C)\longrightarrow \R$ as follows
\[
H(t, \alpha_q, u_q)=\langle \alpha_q ,\Gamma(t, u_q)\rangle-L(t, u_q)
\]
where $\alpha_q\in T^*_q Q$ and $(t, u_q)\in \pi^{-1}(q)$.  Therefore, the coordinate expression of $H$ is (\ref{hamiltonian}).

Let $\omega_Q=-d\theta_Q$ be the canonical symplectic form on $T^*Q$, where $\theta_Q$ is the Liouville form,  and consider the canonical projection $\pi_1:\R\times (T^*Q \times_Q C)\longrightarrow T^*Q$.
Define the 2-form $\Omega_H$ on $\R\times (T^*Q\times_Q C)$ by $\Omega_H=\pi^*_1\omega_Q+dH\wedge dt$. Then, $(dt, \Omega_H)$ is a precosymplectic structure on $\R \times (T^*Q \times_Q C)$ (see \cite{LMM}). 

Eqs. (\ref{eqH}) and (\ref{constraint}) can be intrinsically written as 
\begin{equation}\label{qqq}
i_X\Omega_H=0,\quad i_X dt=1
\end{equation}
Since $(dt, \Omega_H)$ is a precosymplectic structure, Eqs. (\ref{qqq}) need not have a solution, in general. 

Applying the Dirac-Bergmann-Gotay-Nester algorithmn \cite{Dir,Got} to the precosymplectic system $$(\R\times (T^*Q \times_Q C), dt, \Omega_H, H)$$ (see \cite{DLM}) we obtain that Eqs. (\ref{constraint}) correspond to the primary constraints for the precosymplectic system: 
\[
\phi^a=\frac{\partial H}{\partial u^a}=0
\]
Eqs. (\ref{qqq}) have algebraic solution along the first constraint submanifold $P_0$ determined by the vanishing of the primary constraints. On the points of $P_0$ there is at least a pointwise solution   of Eq. (\ref{qqq}), but such solutions are not, in general, tangent to $P_0$. These points must be removed leaving a subset $P_1\subset P_0$ (it is assumed than $P_1$ also is   a submanifold).  Thus, we have to restrict  to a submanifold $P_2$ where the solutions of (\ref{qqq}) are tangent to $P_1$. Proceeding further this way, we obtain a sequence of submanifolds
\[
\cdots\hookrightarrow P_k\hookrightarrow\cdots\hookrightarrow P_2\hookrightarrow P_1\hookrightarrow
P_0\hookrightarrow \R\times (T^*Q\times_Q C)
\]
If this algorithm stabilizes, i.e. there exists a positive integer $k\in \ene$ such that $P_k=P_{k+1}$ and $\dim P_k\not= 0$, then we shall obtain a final submanifold $P_f=P_k$, on which a vector field $X$ exists such that
\begin{equation}\label{final}
( i_X \Omega_H)_{|P_f}=0, \quad (i_X dt=1)_{|P_f}
\end{equation}
The constraints determining $P_f$ are known, in the control literature, as higher order conditions for optimality.

If $X$ is a solution of (\ref{final}) then every arbitrary solution on $P_f$  is of the form $X'=X+ \xi$, where $\xi\in (\ker \Omega_H\cap \ker dt)\cap TP_f$.

Therefore, a necessary condition for optimality of the curve  $\gamma: \R\rightarrow \R\times C$, $\gamma(t)=(t, q(t), u(t))$  is the existence of a lift $\tilde{\gamma}$ of $\gamma$ to $P_f$ such that $\tilde{\gamma}$ is an integral curve of a solution  to Eqs. (\ref{final}). 

In the regular case, the final constraint manifold will be $P_0$ (that is, $P_0=P_f$) and all the constraints are of the second kind following the  classification of Dirac (see \cite{LMM}). In such  case, $(P_0, \Omega, \eta )$ is a cosymplectic manifold, where $\Omega$ and $\eta$ denote the restrictions of $\Omega_H$ and $dt$ to the submanifold $P_0$. 
Denote also by $\omega$ and $\theta$ the restrictions of $\pi^*_1\omega_Q$ and $\pi^*_1\theta_Q$  to $P_0$.

The cosymplecticity of $(P_0, \eta, \Omega)$ is locally equivalent to the regularity of the matrix
\[
\left(
\frac{\partial^2 H}{\partial u^a\partial u^b}
\right)_{1\leq a,b\leq m}
\]  
along  $P_0$.
The dynamical equations for the optimal control problem will become
\begin{equation}\label{aqq}
i_{X}\Omega=0,\quad i_{X}\eta=1
\end{equation}
Taking coordinates $(t, q^A, p_A)$ on $P_0$, then (\ref{aqq}) are equivalent to: 
\begin{equation}\label{eqH1}
\left\{\begin{array}{l}
\dot{q}^A(t)=\displaystyle{\frac{\partial H_{|P_0}}{\partial p_A}(t, q(t), p(t))}\\
$\,$\\
\dot{p}_A(t)=-\displaystyle{\frac{\partial H_{|P_0}}{\partial q^A}(t, q(t), p(t))}\; ,\\
\end{array}\right.
\end{equation}
 where  we have substituted in (\ref{eqH}) the control variables $u^a$ by its value $\bar{u}^a=f^a(t, q, p)$, applying the Implicit Function Theorem to the primary constraints $\phi^a=0$. 
This also implies that we have a canonical projection from $P_0$ onto $\R$, say $\pi_0: P_0\rightarrow \R$ .

In such  case, there exists a unique solution  $X_{P_0}$ of Eq. (\ref{aqq}):
\[
i_{X_{P_0}}\Omega=0,\quad i_{X_{P_0}}\eta=1
\]
 and its flow preserves the cosymplectic structure given by $\Omega$ and $\eta$.
That is, if we denote by $F_{h}$ the flow of $X_{P_0}$ then 
$F_{h}^*\Omega=\Omega$ and $F_{h}^*\eta=\eta$.
In local coordinates, $F_h(t_0, q_0, p_0)=(t_0+h, {q}_1, p_1)$. Denote by $F_h^{(2)}$ the mapping  $F^{(2)}_h(t_0, q_0, p_0)=(q_1, p_1)$, and  by
$F_{t_1, t_0}: P^{t_0}_0\longrightarrow  P^{t_1}_0$ the mapping defined by
\[
F_{t_1, t_0}(q_0, p_0)=F^{(2)}_{t_1-t_0}(t_0, q_0, p_0)\; ,
\]
where we write $P^{t}_0=(\pi_{0})^{-1}(t)$, with $t\in\R$.
Obviously, $F_{t_2, t_1}\circ F_{t_1, t_0}=F_{t_2, t_0}$ in their common domain.

The submanifolds $P^t_0$ naturally inherit a symplectic structure $\omega_t$ by taking the restriction of $\omega$ to $P^t_0$. Similarly, denote  by $\theta_t$ the restriction of $\theta$ to $P^{t}_0$,  then $\omega_{t}=-d\theta_t$.

It is easy to deduce that, in such  case, $F_{t_1, t_0}$ is a symplectomorphism; that is, $F_{t_1, t_0}^*\omega_{t_1}=\omega_{t_0}$, noting  that 
\[
\Omega=\omega+dH_{|P_0}\wedge \eta
\]
This last remark will be interesting for constructing geometrical integrators for explicitly time-dependent optimal control systems.

\section{Generating functions}\label{section}

Let $(M_{i}, \omega_{i})$, $i=0,1$ be two exact  symplectic manifolds (i.e. $\omega_i$ is symplectic and exact, $\omega_{i}=-d\theta_{i}$, $i=0, 1$)  and suppose that $g: M_{0}\rightarrow M_{1}$ is a diffeomorphism. Denote by $\hbox{Graph}(g)$ the graph of $g$, $\hbox{Graph}(g)=\{(x_0, g(x_0))\; /\; x_0\in M_0\} \subset M_0\times M_1$.
Denote by $\pi_i: M_0\times M_1\rightarrow M_i$, $i=0,1$ the canonical projections, and  consider the 1-form and 2-form on $M_0\times M_1$ defined by
\begin{eqnarray*}
\Theta_{(1,0)}&=&\pi_1^*\theta_1-\pi_0^* \theta_0\\
{\Omega}_{(1,0)}&=&\pi_1^*\omega_1-\pi_0^* \omega_0=-d{\Theta}_{(1,0)}\\
\end{eqnarray*}
As it is well known $\Omega_{(1,0)}$ is a symplectic form.   

Let  $i_g: \hbox{Graph}(g)\hookrightarrow M_0\times M_1$ be  the inclusion map, then
\[
i_g^*{\Omega}_{(1,0)}=({\pi_0}_{|\hbox{\small Graph}(g)})^*(g^*\omega_1-\omega_0)
\]
Using this equality, it is clear that  $g$ is a symplectomorphism  if and only if $i_g^*{\Omega}_{(1,0)}=0$, that is, if $\hbox{Graph}(g)$ is a Lagrangian submanifold of $(M_0\times M_1, {\Omega}_{(1,0)})$.

Now, if $g$ is a symplectomorphism we have
\[
i_g^*{\Omega}_{(1,0)}=-d i_g^*\Theta_{(1,0)}=0
\]
and, therefore, at least locally, there exists a function $S: \hbox{Graph}\, (g) \rightarrow \R$ such that
\begin{equation}\label{eee}
i_g^*\Theta_{(1,0)}=dS
\end{equation}
Let $(q_0, p_0)$ and $(q_1, p_1)$ Darboux coordinates in $M_0$ and $M_1$, respectively. Since $\hbox{Graph}(g)$ is diffeomorphic to $M_0$, we can take   $(q_0,p_0)$ as natural coordinates in $\hbox{Graph}(g)$.  Since  $(q_0, p_0, q_1, p_1)$ are coordinates in $M_0\times M_1$, then,  along $\hbox{Graph}(g)$, we have $q_1=q_1(q_0,p_0)$, $p_1=p_1(q_0,p_0)$ and
\[
p_1\,dq_1-p_0dq_0=dS(q_0,p_0) 
\]

\subsection{Generating functions of the first kind}

Assume that in a neighborhood of some point $x\in \hbox{Graph}(g)$, we can change this system of coordinates by new independent coordinates 
$(q_0, q_1)$ (the local condition is that $\det \left( \partial q_1/\partial p_0\right)\not=0$).
 In such  case, the function $S$ can be  expressed locally  as 
$
S=S(q_0,p_0)=S_1(q_0,q_1)$.

\begin{definition}
The function $S_1(q_0,q_1)$ will be called a {\bf generating function of the first kind} of the symplectomorphism $g$.
\end{definition}
From (\ref{eee}) we deduce that
\begin{equation}\label{aq10}
\left\{
\begin{array}{l}
\displaystyle{p_0= -\frac{\partial S_1}{\partial q_0}}\\
\vspace{-0.3cm}$\,$\\
\displaystyle{p_1= \frac{\partial S_1}{\partial q_1}}\\
\end{array}
\right.
\end{equation}
(see Fig. 1).

Conversely, if $S_1(q_0, q_1)$ is a function such that 
$
\det \left(\frac{\partial^2 S_1}{\partial q_0\partial q_1}\right)\not= 0
$
then $S_1(q_0, q_1)=(p_0,p_1)$ is a generating function of some canonical transformation $g$ implicitly determined by Eqs. (\ref{aq10}), $g(q_0, p_0)=(q_1, p_1)$ (see \cite{Arno}).

\

\unitlength=0.8mm
\special{em:linewidth 0.4pt}
\linethickness{0.4pt}
\begin{picture}(80.33,50.00)(-40,0)
\put(20.00,25.00){\makebox(0,0)[cc]{$(q_0,q_1)$}}
\put(27.00,25.00){\vector(1,0){38.00}}
\put(68.00,25.00){\makebox(0,0)[lc]{$(q_0, -\frac{\partial S_1}{\partial q_0}, q_1, \frac{\partial S_1}{\partial q_1})$}}
\put(80.00,0.00){\makebox(0,0)[cc]{$(q_1, \frac{\partial S_1}{\partial q_1})$}}
\put(80.00,50.00){\makebox(0,0)[cc]{$(q_0, -\frac{\partial S_1}{\partial q_0})$}}
\put(80.00,30.33){\vector(0,1){16.67}}
\put(80.00,20.00){\vector(0,-1){16.67}}
\put(26.33,28.00){\vector(2,1){43.33}}
\put(26.00,21.67){\vector(2,-1){43.00}}
\put(81.33,37.67){\makebox(0,0)[lc]{$(\pi_1)_{|\hbox{\small Graph}(g)}$}}
\put(81.33,11.67){\makebox(0,0)[lc]{$(\pi_2)_{|\hbox{\small Graph}(g)}$}}
\put(50.00, -5){\makebox(0,0)[cc]{Fig.1}}

\end{picture}

\

Now suppose  that $M$ is a fiber bundle over the real line $\R$, $\pi: M\rightarrow \R$, and $M_t=\pi^{-1}(t)$ are the fibers, where each fiber  $M_t$ is equipped with a symplectic form $\omega_t$.
Let $g_{(s, t)}: M_t\rightarrow M_s$
 be a  two-parameter family of symplectomorphisms 
satisfying 
\[
g_{(t_2, t_1)}\circ g_{(t_1, t_0)}=g_{(t_2, t_0)}
\]
Next, we shall show how this composition law can be  translated in terms of  their respective generating functions. Moreover, 
the following results will give a geometric interpretation of the Discrete Euler-Lagrange equations \cite{Mars6}.

\begin{theorem}\label{teo1}
Let  $S_1^{(t_N, t_0)}$ be a function defined by
\[
S_1^{(t_N, t_0)}(q_0, q_{N})=\sum_{k=0}^{N-1}S_1^{(t_{k+1}, t_k)}(q_k, q_{k+1})
\]
where $q_k\in M_{t_k}$, $1\leq k\leq N-1$, are stationary points of the right-hand side, that is
\begin{eqnarray*}
0&=&D_2 S^{(t_k, t_{k-1})}_1(q_{k-1}, q_k)+D_1 S^{(t_{k+1}, t_k)}_1(q_k, q_{k+1}),
 \quad 1\leq k\leq N-1. 
\end{eqnarray*}
If $S^{(t_k, t_{k-1})}_1$ are generating functions of the first kind for $g_{(t_{k}, t_{k-1})}$, then  $S^{t_N,t_0)}_1$ is a generating function of the first kind for $g_{(t_N, t_0)}: M_{t_0}\rightarrow M_{t_N}$.
\end{theorem}
{\bf Proof:}
Recursively, it is suffices to give the proof for $N=2$:
\[
S_1^{(t_2,t_0)}(q_0, q_2)=S^{(t_1,t_0)}_1(q_0, x)+S^{(t_2, t_1)}_1(x, q_2)
\]
where $x$ is an stationary point of the right-hand side.

From the definitions of  generating functions for $g_{(t_2,t_1)}$ and $g_{(t_1, t_0)}$ 
\begin{eqnarray*}
p_1\,dq_1-p_0\, dq_0&=&dS^{(t_1, t_0)}(q_0, q_1)\\
p_{2}\,dq_{2}-p_1\, dq_1&=&dS^{(t_2, t_1)}_1(q_1, q_2)\\
\end{eqnarray*}
and therefore 
\[
p_2\,dq_2-p_0\, dq_0=d(S^{(t_2, t_1)}_1(q_0, q_1)+S^{(t_1, t_0)}_1(q_1, q_2))
\]
It follows that
\[
0=D_2 S^{(t_1, t_0)}_1(q_{0}, q_1)+D_1 S^{(t_2, t_1)}_1(q_1, q_{2})
\]
and, obviously, for this choice of $q_1$ then 
\[
S^h_1(q_0, q_1)+S^h_1(q_1, q_2)
\]
is a generating function of the first kind of $g_{(t_2, t_0)}$.\bull

Now, we are in condition to bring this procedure to the limit when the number of 
subintervals increases to infinity.  Consider as its continuous counterpart a cosymplectic manifold $(M, \eta, \omega)$, where $M$ is still a fiber bundle over $\R$ ($\pi_{\R}: M\rightarrow \R$) and $\eta=\pi_{\R}^*(dt)$. Denote by $M_t=\pi_{\R}^{-1}(t)$, $t\in \R$. Take  a Hamiltonian function $H: M\rightarrow \R$ and its  Hamiltonian vector field $X_H$
given by
\[
i_{X_H}\omega=0\qquad \hbox{and}\qquad i_{X_H}\eta=1
\]
Let  $F_{(t, s)}: M_{s}\rightarrow M_t$ be  the two-parameter family of symplectomorphisms  generated by $X_H$ (see section \ref{optimal}) and consider as symplectic form on each fiber $M_t$ the restriction of $\omega$ to this fiber. 

We shall give a characterization of the generating functions of the first kind associated to $F_{(t,s)}$ for $t$ close enough to $s$. For doing that, consider Darboux coordinates $(t, q^A, p_A)$ on $M$ and assume the regularity condition 
$\displaystyle{
\det \left(
\frac{\partial^2 H}{\partial p_A\partial p_B}
\right)\not=0
}
$. Thus, 
\begin{proposition}\label{proposition1}
A generating function of the first kind for $F_{(t, s)}$ is given by
\[
S_1^{(t_1, t_0)}(q_0,q_1)=\int^{t_1}_{t_0} \left(p(t)\dot{q}(t)-H(t, q(t), p(t))\right) \, dt
\]
where $t\rightarrow (t, q(t), p(t))$ is an integral curve of the Hamilton equations such that $q(t_0)=q_0$ and $q(t_1)=q_1$.
\end{proposition}

{\bf Proof:}
We only use  Hamilton equations and integration by parts:
\begin{eqnarray*}
\frac{\partial S^{(t_1, t_0)}_1}{\partial q_0}(q_0, q_1)&=&
\int^{t_1}_{t_0}\left(\frac{\partial p}{\partial q_0}\dot{q}+p\frac{\partial 
\dot{q}}{\partial q_0}-\frac{\partial H}{\partial q}\frac{\partial q}{\partial q_0}-\frac{\partial H}{\partial p}\frac{\partial p}{\partial q_0}\right)\, dt\\
&=&%
\int^{t_1}_{t_0}\left({p}\frac{\partial \dot q}{\partial q_0}+\dot{p}\frac{\partial q}{\partial q_0}\right)\,dt
\\
&&=-p_0+p_1\frac{\partial q_1}{\partial q_0}=-p_0\\
\end{eqnarray*}
and 
\begin{eqnarray*}
\frac{\partial S^{(t_1,t_0)}_1}{\partial q_1}(q_0, q_1)&=&
\int^{t_1}_{t_0}\left(\frac{\partial p}{\partial q_1}\dot{q}+p\frac{\partial 
\dot{q}}{\partial q_1}-\frac{\partial H}{\partial q}\frac{\partial q}{\partial q_1}-\frac{\partial H}{\partial p}\frac{\partial p}{\partial q_1}\right)\, dt\\
&=&%
\int^{t_1}_{t_0}\left(p\frac{\partial q}{\partial q_1}+\dot{p}\frac{\partial q}{\partial q_1}\right)\,dt
\\
&&=p_1-p_0\frac{\partial q_0}{\partial q_1}=p_1\bull\\
\end{eqnarray*}

\begin{remark}
{\rm 
Suppose that $t_{i+1}-t_i=h$, for all $i=0,\cdots N-1$, then from Theorem \ref{teo1} we have 
\[
S^{Nh}_1(q_0, q_N)=\sum_{k=0}^{N-1}S_1^{h}(q_k, q_{k+1})
\]
where
\begin{eqnarray*}
0=D_2 S^{h}_1(q_{k-1}, q_k)+D_1 S^{h}_1(q_k, q_{k+1}),
 \quad 1\leq k\leq N-1. 
\end{eqnarray*}
Now, if we take as new generating function an adequate approximation $S^h_d$ of $S_1^h$ then 
\begin{eqnarray*}
0=D_2 S^{h}_d(q_{k-1}, q_k)+D_1 S^{h}_d(q_k, q_{k+1}),
 \quad 1\leq k\leq N-1. 
\end{eqnarray*}
are the well-known Discrete Euler-Lagrange equations (see \cite{Mars6} and references therein). For instance, one can take
\[
S_d^h(q_0, q_1)=h {\cal L}(\alpha q_0+(1-\alpha) q_1, \frac{q_1-q_0}{h}),\quad 
\alpha\in [0,1]
\]
or alternatively, we could have considered more accurate approximations. Here, we are assuming that ${\cal L}: \R\times TQ\rightarrow\R$ is a  Lagrangian function related via Legendre transformation with the Hamiltonian function $H$ (see \cite{Arno}) which is locally possible because of the regularity of $H$.
}
\end{remark}

Denote by $S_1(q_0, q_1, t_0, t_1)=S_1^{(t_1, t_0)}(q_0, q_1)$. From Proposition (\ref{proposition1}), it is easy to show that:
\begin{eqnarray*}
D_3 S_1(q_0, q_1, t_0, t_1)=D_3 S^{(t_1, t_0)}(q_0, q_1)&=&H(t_0, q_0, p_0)\\
D_4 S_1(q_0, q_1, t_0, t_1)=D_4 S^{(t_1, t_0)}(q_0, q_1)&=&-H(t_1, q_1, p_1)
\end{eqnarray*}
(see also \cite{Mars6}).
As a consequence
\begin{equation}\label{pou}
D_4 S^{(t_k, t_{k-1})}(q_{k-1}, q_k)+D_3 S^{(t_{k+1}, t_k)}(q_k, q_{k+1})=0
\end{equation}

It should be noticed that if we take a new  function $S_d^{(t_{k+1}, t_k)}$ as an adequate approximation  of $S^{(t_{k+1}, t_k)}$, 
then solutions $\{q_0, q_1, \ldots, q_N\}$ of equations
\begin{eqnarray*}
D_2 S^{(t_k, t_{k-1})}_d(q_{k-1}, q_k)+D_1 S^{(t_{k+1}, t_k)}_d(q_k, q_{k+1})=0,
 \quad 1\leq k\leq N-1. 
\end{eqnarray*}
 do not satisfy  (\ref{pou}) for arbitrary values of $t_{k-1}, t_k, t_{k+1}$.
Therefore, we may write the system of difference equations
\begin{equation}\label{pouy}
\left\{
\begin{array}{l}
D_2 S^{(t_k, t_{k-1})}_d(q_{k-1}, q_k)+D_1 S^{(t_{k+1}, t_k)}_d(q_k, q_{k+1})=0,\\
D_4 S^{(t_k, t_{k-1})}_d(q_{k-1}, q_k)+D_3 S^{(t_{k+1}, t_k)}_d(q_k, q_{k+1})=0,
\end{array}
\right.
\end{equation}
which under regularity assumptions will determine a time-dependent discrete flow
$$\Phi(q_{k-1}, q_k, t_{k-1}, t_k)=(q_k, q_{k+1}, t_k, t_{k+1})$$
with variable step size $h_k=t_{k+1}-t_k$ (see \cite{Kane1,Lee1,Lee2,LD1,Mars6}).

\subsection{Generating functions of the second kind}

The construction of more general generating functions will be useful in  next sections. For instance,  suppose that $(q_0, p_1)$ are independent local coordinates on $\hbox{Graph}(g)$. Then the function $S$ is written as $S=S(q_0,p_1)$.

We have
\[
p_1\,dq_1-p_0\, dq_0=-q_1\, dp_1 + d(q_1 p_1)-p_0\, dq_0=dS.
\]
If we define 
\[
S_2(q_0,p_1)=q_1 p_1-S(q_0,p_1),
\]
where $q_1$ is expressed in terms of $q_0$ and $p_1$, then we deduce that 
\[
q_1\, d p_1+p_0 dq_0=dS_2(q_0, p_1)
\]

 \begin{definition}
The function $S_2(q_0,p_1)$ will be called a {\bf generating function of the second kind} of the symplectomorphism $g$.
\end{definition}

We have that
\begin{equation}\label{aq1}
\left\{
\begin{array}{l}
\displaystyle{p_0= \frac{\partial S_2}{\partial q_0}}\\
\vspace{-0.3cm}$\,$\\
\displaystyle{q_1= \frac{\partial S_2}{\partial p_1}}\\
\end{array}
\right.
\end{equation}

Conversely, if $S_2(q_0,p_1)$ is a generating function such that 
$
\displaystyle{\det \left(\frac{\partial^2 S_2}{\partial q_0\partial p_1}\right)\not= 0}
$
then $S_2$ is a generating function of some local symplectomorphism determined by Eqs. (\ref{aq1}) (see \cite{Arno}).

Denote by $F_{(t_2, t_1)}: M_{t_1}\rightarrow M_{t_2}$ the two-parametric group of canonical transformations generated by the Hamiltonian vector field $X_H$, as in the preliminaries to Proposition \ref{proposition1}. We have the following. 

\begin{theorem}\label{teo2}
Let a function $S_2^{(t_N, t_0)}$ be defined by
\begin{equation}\label{mo3}
S_2^{(t_N,t_0)}(q_0, p_{N})=\sum_{k=0}^{N-1}S^{(t_{k+1}, t_k)}_2(q_k, p_{k+1})-\sum_{k=1}^{N-1}q_{k}p_{k}
\end{equation}
where $q_k$, $1\leq k\leq N$, and $p_k$, $0\leq k\leq N-1$, are stationary points of the right-hand side, that is
\begin{eqnarray}
q_k&=&\frac{\partial S^{(t_{k-1}, t_k)}_2}{\partial p}(q_{k-1}, p_k), \quad 1\leq k\leq N,\label{mo1}\\ 
p_k&=&\frac{\partial S^{(t_{k}, t_{k+1})}_2}{\partial q}(q_{k}, p_{k+1}), \quad 0\leq k\leq N-1,\label{mo2} 
\end{eqnarray}
then $S^{(t_N, t_0)}_2$ is a generating function of the second kind for $F_{(t_N, t_0)}: M_{t_0}\rightarrow M_{t_N}$.
\end{theorem}
{\bf Proof:}
It follows as in Theorem \ref{teo1}.\bull

As a consequence, we have that 
\begin{equation}\label{mo8}
S^{(t_N, t_0)}(q_0,p_N)=q_Np_N-S^{(t_N, t_0)}_2(q_0,p_N)=\sum_{k=0}^{N-1}\left[q_{k+1}p_{k+1}-S^{(t_{k+1}, t_k)}_2(q_k, p_{k+1})\right]\; ,
\end{equation}
where the unknown coordinates are given by (\ref{mo1}) and (\ref{mo2}).

\begin{proposition}\label{proposition2}
A generating function of the second kind for $F_{(t_1, t_0)}$ is given by
\[
S_2^{(t_1, t_0)}(q_0,p_1)=p_1q_1-\int^{t_1}_{t_0} \left(p(t)\dot{q}(t)-H(t, q(t), p(t))\right)\, dt
\]
where $t\rightarrow (q(t), p(t))$ is an integral curve of the Hamilton equations such that $q(t_0)=q_0$ and $p(t_1)=p_1$.
\end{proposition}

{\bf Proof:}
It is proved in a similar way to Proposition \ref{proposition1}.\bull


Denote by $S_2(t, q_0, p_1)=S_2^{(0,t)}(q_0, p_1)$ then it is easy to show that (see, for instance  \cite{Hair})

\begin{theorem}[Hamilton-Jacobi equation for $S_2$]
If $S_2(t, q_0, p_1)$ is a solution of the partial differential equation
\begin{equation}\label{azz}
\frac{\partial S_2}{\partial t}=H(\frac{\partial S_2}{\partial p_1}(t, q_0, p_1), p_1),\qquad S_2(0, q_0, p_1)=q_0 p_1
\end{equation}
then the mapping $(q_0, p_0)\longrightarrow (q_1, p_1)$ defined by Eqs. (\ref{aq1}) is the exact flow of the Hamiltonian 
system determined by $H$. 
\end{theorem}

\section{Optimal control of Discrete-time systems}

In this section we shall define the general solution  of an optimization problem for discrete systems and analyze its geometric behaviour, in particular, the symplecticity.

Suppose that the discrete state equations are given by the dynamical equation
\begin{equation}\label{eq1}
q^A_{k+1}=f^A(k,q_k, u_k), \quad k=0, 1, \ldots, N-1, \quad A=1, 2,\ldots, m
\end{equation}
or, shortly,  
$
q_{k+1}=f(k,q_k, u_k)$, 
where $q_0$ is initially given.

The associate performance index or objective function is:
\begin{equation}\label{eq2}
J=\bar{S}(N, q(N))+\sum_{k=0}^{N-1} \bar{L}(k, q_k, u_k)
\end{equation}
where $\bar{S}$ is a function of the final time and   state at the final time $N$, and $\bar{L}$ is time-varying function of the state and control input at each intermediate discrete time $k$.

The optimal control problem is solved finding  controls $u^*_k$, $k=0, 1,\ldots N-1$, that drive the system along a trajectory $q_k^*$, $k=0, 1,\ldots, N$, verifying the state equations such that the performance index is minimized.

\subsection{Problem solution}

Let us now solve the optimal control problem for the discrete optimal problem determined by (\ref{eq1}) and (\ref{eq2}) using the Lagrange multiplier approach.
Considering the state Eqs. (\ref{eq1}) as constraint equations,  then we have $N\cdot m$ constraints, and we associate  a Lagrange multiplier to each constraint.
Next, we  construct the augmented performance index $J'$ by
\begin{equation}\label{eq3}
J'=\sum_{k=0}^{N-1}\left[ p_{k+1}(f(k, q_k, u_k)-q_{k+1})-\bar{L}(k, q_k, u_k)\right]-\bar{S}(N, q(N))
\end{equation}
where $p_{k+1}=((p_{k+1})_A)$ are considered as Lagrange multipliers with $A=1,\ldots, n$ and $k=0, \ldots, N-1$.

Taking the Hamiltonian function
\[
\bar{H}(k, q_k, p_{k+1}, u_k)=p_{k+1}f(k, q_k, u_k)-\bar{L}(k, q_k,u_k)
\] 
we deduce that the  necessary conditions for a constrained minimum are thus given by:
\begin{eqnarray}
q_{k+1}&=& \frac{\partial \bar{H}}{\partial p}(k, q_k, p_{k+1}, u_k)=f(k, q_k, u_k) \label{er1}\\ 
p_{k}&=& \frac{\partial \bar{H}}{\partial q}(k, q_k, p_{k+1}, u_k)=p_{k+1}\frac{\partial f}{\partial q}(k, q_k, u_k)-\frac{\partial \bar{L}}{\partial q}(k, q_k, u_k)\label{er2}\\ 
0&=&\frac{\partial \bar{H}}{\partial u}(k, q_k, p_{k+1}, u_k)=p_{k+1}\frac{\partial f}{\partial u}(k, q_k, u_k)-\frac{\partial \bar{L}}{\partial u}(k, q_k, u_k)\label{er3}
\end{eqnarray}
where $0\leq k\leq N-1$,
and the transversality conditions 
\[
p_N=-\frac{\partial \bar{S}}{\partial q}(N, q_N) \quad \hbox{and}\quad q_0\quad \hbox{fixed}.
\]

Observe that the recursion for the state $q_k$ develops forward in time, but the co-state variable $p_k$ develops backwards in time. Therefore the required boundary conditions for finding a solution are the initial  state $q_0$ and the final co-state $p_N$. 

Assume that 
\[
\det \left(\frac{\partial^2 \bar{H}}{\partial u_a \partial u_b}\right)\not=0
\]
then, locally, $u^*_k=h(k, q_k, p_{k+1})$.  
If we denote, by 
\[
\tilde{H}(k, q_k, p_{k+1})=\bar{H}(k, q_k, p_{k+1}, u^*_k)
\]
then Eqs. (\ref{er1}), (\ref{er2}) are rewritten as 
\begin{eqnarray}
q_{k+1}&=& \frac{\partial \tilde{H}}{\partial p}(k, q_k, p_{k+1}) \label{err1}\\ 
p_{k}&=& \frac{\partial \tilde{H}}{\partial q}(k, q_k, p_{k+1})\label{err2} 
\end{eqnarray}
with $0\leq k\leq N_1$.

Consider the function
\begin{eqnarray*}
G_k(q_k, q_{k+1}, p_{k+1})&=&\tilde{H}(k, q_k, p_{k+1})-p_{k+1}q_{k+1}, \quad 0\leq k\leq N-1.\\
\end{eqnarray*}
Then, for a fixed $k$: 
\begin{eqnarray*}
dG_k&=&\frac{\partial \bar{H}}{\partial q_k}(k, q_k, p_{k+1})\,
dq_k+\frac{\partial \bar{H}}{\partial p_{k+1}}(k, q_k, p_{k+1})\,
dp_{k+1}
-p_{k+1}\, dq_{k+1}-q_{k+1}\, dp_{k+1}\; .\\
\end{eqnarray*}
Along  solutions of Eqs. (\ref{er1}), (\ref{er2}) and (\ref{er3}) we have:
\begin{eqnarray*}
{dG_k}_{|\hbox{\footnotesize Sol}}&=& p_k\, dq_k-p_{k+1}\, d q_{k+1}\; ,
\end{eqnarray*}
which implies
\begin{equation}\label{symp}
dp_{k}\wedge dq_k=dp_{k+1}\wedge dq_{k+1}\; .
\end{equation}
along the solution of (\ref{er1})-(\ref{er3}).

In the next subsection, we shall analyze the geometric meaning   of Eq. (\ref{symp}), which it is obviously interpreted as symplecticity of discrete optimal control problems in terms of a natural symplectic form.

\subsection{Generating functions of the second kind and discrete optimal control problems}

From Proposition \ref{proposition1} the following function is a generating function of the second kind for the cosymplectic Hamiltonian system 
$(P_0, \eta, \Omega, H_{|P_0})$, which determines the dynamics of the optimal control problem given by (\ref{general}) and (\ref{general1}): 
\begin{equation}\label{ppp}
S^{(t_1, t_0)}_2(q_0, p_1)=p_1q_1-\int_{t_0}^{t_1} \left(p(t)\dot{q}(t)-H_{|P_0}(t, q(t),p(t))\right)\,dt\; ,
\end{equation}
where $t\rightarrow (t, q(t), p(t))$ is the integral curve on $P_0$ of the vector field $X_{P_0}$. Here $X_{P_0}$ is the unique solution of equation
\[
i_{X_{P_0}}\Omega=dH_{|P_0}\qquad i_{X_{P_0}}\eta=1
\]
with $(q(t_0), p(t_0))=(q_0, p_0)$ and $(q(t_1), p(t_1))=(q_1, p_1)$.

We now focus on the construction of a numerical  integrator for the Hamiltonian system $(P_0, \eta , \Omega, H_{|P_0})$ 
by using an approximation of the generating function. As we shall show, the obtained method  also realize the integration steps by symplectomorphism transformations; then, it is a symplectic integrator.

First take  a fixed time interval $h=t_{k+1}-t_k$, $k=0, \ldots, N-1$. 

Assume that we are working  on vector spaces, and consider the following natural approximation:
\begin{eqnarray*}
\tilde{S}^{h}_2(k, q_k, p_{k+1})&=&p_{k+1}q_{k+1}-hp_{k+1}\left(\frac{q_{k+1}-q_k}{h}\right)-h\tilde{L}(k, q_k, p_{k+1})\\
&&+hp_{k+1}\tilde{\Gamma}(k, q_k, p_{k+1})
\end{eqnarray*}
where, for instance,  $\tilde{L}(k, q_k, p_{k+1})=L_{|P_0}(t_0+kh, q_k, p_{k+1})$ and 
$\tilde{\Gamma}(k, q_k, p_{k+1})=\Gamma_{|P_0}(t_0+kh, q_k, p_{k+1})$.

If we denote by $\tilde{f}(k, q_k, p_{k+1})$ the function
\begin{equation}\label{ase}
\tilde{f}(k, q_k, p_{k+1})=h\tilde\Gamma(k, q_k, p_{k+1})+q_{k}
\end {equation}
then, 
\[
\tilde{S}^{h}_2(k, q_k, p_{k+1})=p_{k+1}\tilde{f}(k, q_k, p_{k+1})-\tilde{L}(k, q_k, p_{k+1})= \tilde{H}(k, q_k, p_{k+1})\; .
\]
Thus,  equations 
\begin{equation}\label{aq4}
\left\{
\begin{array}{l}
\displaystyle{p_k= \frac{\partial \tilde{S}^{h}_2}{\partial q^k}(k, q_k, p_{k+1})=\frac{\partial \tilde{H}}{\partial q^k}(k, q_k, p_{k+1})}\\
$\,$\\
\displaystyle{q_{k+1}= \frac{\partial \tilde{S}^{h}_2}{\partial p_{k+1}}(k, q_k, p_{k+1})=\frac{\partial \tilde{H}}{\partial p_{k+1}}(k, q_k, p_{k+1})}
\end{array}
\right.
\end{equation}
are exactly (\ref{err1}) and (\ref{err2}) and the symplecticity condition 
 (\ref{symp}) for discrete optimal control problems is now a trivial consequence of the generating function construction.

\begin{remark}{
\rm
It is also possible to construct symplectic numerical methods of higher order;  for instance, considering better approximations of the Hamilton Jacobi equation (\ref{azz}) (see \cite{Chan} and references therein). Assume for simplicity that the Hamiltonian is autonomous, that is, $H\equiv H(q, p)$. Now, first  expands the generating function $S_2(t, q_0, p_1)$ as:
\[
S_2(t, q_0, p_1)=q_0  p_1+\sum_{i=1}^{\infty}t^iG_i(q_0, p_1),
\]
inserts expression into Hamilton-Jacobi equation (\ref{azz}) and compares equal powers of $t$. This yields
\begin{eqnarray*}
G_1(q_0, p_1)&=&H(q_0, p_1)\\
G_2(q_0, p_1)&=&\frac{1}{2}\left(\frac{\partial H}{\partial q^A_0} \frac{\partial H}{\partial p_{1A}}\right)\\
G_3(q_0, p_1)&=&\frac{1}{6}\left(\frac{\partial^2 H}{\partial p_{1A}\partial p_{1B}}\frac{\partial H}{\partial q_0^A}\frac{\partial H}{\partial q_0^B}+\frac{\partial^2 H}{\partial p_{1A}\partial q_0^B}\frac{\partial H}{\partial q_0^A}\frac{\partial H}{\partial p_{1B}}+\frac{\partial^2 H}{\partial q^A_0\partial q^B_0}\frac{\partial H}{\partial p_{1A}}\frac{\partial H}{\partial p_{1B}}\right)\\
\cdots&=&\cdots
\end{eqnarray*}
Using the truncated series, we obtain an approximated generating function:
\[
S_2^h(q_k, p_{k+1})=q_k\cdot p_{k+1}+\sum_{i=1}^{r}h^rG_i(q_k, p_{k+1})
\]
which defines a symplectic method of order $r$. 

Other approaches are also admissible without using higher derivatives of the Hamiltonian $H$, for instance, symplectic  or symplectic partitioned Runge-Kutta methods (see \cite{Hair,Sanz}).
}
\end{remark}

\section{Discrete Hamiltonian systems}

In \cite{EY} Erbe and Yan have considered discrete linear Hamiltonian systems of the form: 
\begin{eqnarray*}
\Delta y(t)&=& B(t) y(t+1) + C(t) z(t)\\
\Delta z(t)&=& -A(t) y(t+1)-B^T (t) z(t)\\
\end{eqnarray*}
where $A, C$ are symmetric and $I-B$ is invertible. Here $\Delta y(t)=y(t+1)-y(t)$, $\Delta z(t)=z(t+1)-z(t)$ and $y, z\in\R^d$.

This problem is a particular case of a discrete Hamiltonian systems of the form
\begin{eqnarray}
\Delta y(t)&=& H_z (t, y(t+1), z(t))\label{fg}\\
\Delta z(t)&=& -H_y (t, y(t+1), z(t))\label{fg1}
\end{eqnarray}
where $\displaystyle{H(t, y, z)=\frac{1}{2}(y^T, z^T)
\left(
\begin{array}{cc}
A(t)&B^T(t)\\
-B(t)&C(t)
\end{array}
\right)
\left(
\begin{array}{c}
y\\
z
\end{array}
\right)
}
$.
The symplecticity of the discrete linear Hamiltonian system was fully studied (see \cite{EY}, for instance, and references therein). The existence of a corresponding symplectic structure for discrete nonlinear Hamiltonian systems given by (\ref{fg}) and (\ref{fg1}) was proposed by Ahlbrandt as an open problem (\cite{Ahl} and also \cite{Shi}).

From the point of view of section \ref{section}, this open problem is easily solved considering as generating function of the second kind the following one: 
\[
S^{(t+1, t)}_2(y(t+1), z(t))=z(t) y(t+1)- H(t, y(t+1), z(t))\; .
\]
Then Eqs. (\ref{fg}) and (\ref{fg1}) are precisely
\[
\left\{
\begin{array}{rcl}
y(t)&=&\frac{\partial S^{(t+1,t})}{\partial z}(y(t+1), z(t))\\
z(t+1)&=& \frac{\partial S^{(t+1,t)}}{\partial y}(y(t+1), z(t))\; ,
\end{array}
\right.
\]
which guarantees the symplecticity of the discrete Hamiltonian system. In order to find the canonical transformation associated to this generating map it is only necessary to impose the local condition (see \cite{Arno}):
  \[
\det \left( \frac{\partial^2 S^{(t+1,t)}_2(y(t+1), z(t)}{\partial y\partial z}\right)\not=0
\]
Then, in a neighbourhood of a point satisfying the above condition, there exists a  symplectomorphism defined by Eqs. (\ref{fg}) and (\ref{fg1}).   

\section*{Acknowledgments}

This work has been  supported by grant BFM2001-2272 (Ministry of Science and Technology, Spain). A. Santamar{\'\i}a Merino wishes to thank the Programa de Formaci\'on de Investigadores of the Departamento de Educaci\'on, Universidades e Investigaci\'on of the Basque Government (Spain) for financial support.


\begin{thebibliography}{99}


\bibitem[Ahlb:93]{Ahl} Ahlbrandt C D 1993 Equivalence of Discrete Euler Equations and Discrete Hamiltonian Systems,   
{\it J. Math. Anal. Appl.} {\bf 180}, 498-478

\bibitem[Arn:78]{Arno} Arnold V I 1978 {\it Mathematical Methods of Classical Mechanics} (Graduate Text in Mathematics 60, Springer-Verlag New York)

 
\bibitem[BobSus:99a]{Bobe} Bobenko A I  and Y B Suris 1999 Discrete Lagrangian reduction, discrete
Euler-Poincar\'e equations, and semidirect products {\it Lett. Math. Phys.} {\bf 49}, 79-93

\bibitem[BobSus:99b]{Bobe1} Bobenko A I and  Suris Y B 1999 Discrete time Lagrangian mechanics on Lie
groups, with an application to the Lagrange top {\it Comm. Math. Phys.} {\bf 204} 
147-188

\bibitem[Cadz:70]{Cadz}  Cadzow J A 1970 Discrete calculus of variations {\it Intern. J. Control.} {\bf 11},
393-407

\bibitem[ChaSco:90]{Chan}  Channell P J and  Scovel C 1990 Symplectic integration of Hamiltonian Systems, {\it Nonlinearity} {\bf 3}, 231-259

\bibitem[ChLeMa:94]{DLM}   Chinea D,  de Le\'on M and   Marrero J C 1994
 The
constraint algorithm for time-dependent Lagrangians,
 {\it J. Math.
Phys.} {\bf 35} (7), 3410-3447



\bibitem[Cort:02]{Cort} Cort\'es J 2002 {\em Geometric, control and numerical aspects of nonholonomic systems} (Lecture 
     Notes in Mathematics, vol. 1793, Springer-Verlag) 


\bibitem[CorMar:01]{JS} Cort\'es J and  Mart{\'\i}nez S 2001  Nonholonomic integrators {\it  
     Nonlinearity} {\bf 14}, 1365-1392

\bibitem[Dirac:64]{Dir}  Dirac P A M 1964  {\it Lecture on Quantum Mechanics} 
(Belfer Graduate School of Science, Yeshiva University, New York)

\bibitem[ErbYan:92]{EY} Erbe L H and Yan P 1992 Disconjugancy for linear Hamiltonian difference systems, {\it J. Math. Anal. Appl.} {\bf 167}, 355-367

\bibitem[GotNes:79]{Got}   Gotay M J and  Nester J M 1979
 Presymplectic
Lagrangian systems I: The Constraint Algorithm and the
Equivalence Theorem,
 {\it Ann. Inst. Henri Poincar\'e}, {\bf A30} 129-142



\bibitem[HaLuWa:02]{Hair}  Hairer E Lubich C and Wanner G 2002 {\it Geometric Numerical  Integration, Structure-Preserving Algorithms for Ordinary Differential Equations} (Springer Series in Computational Mathematics {\bf 31}, Springer-Verlag Berlin  Heidelberg)

\bibitem[JarNor:97a]{Jaro1}  Jaroszkiewicz G and Norton K 1997 Principes of discrete time mechanics I:
Particle systems {\it J. Phys. A}  {\bf 30}, 3115-3144

\bibitem[JarNor:97b]{Jaro2} Jaroszkiewicz G and  Norton K 1997  Principes of discrete time mechanics II:
Classical field theory {\it J. Phys. A } {\bf 30}, 3145-3163

\bibitem[JorPol:64]{Jordan}  Jordan B W and Polak E 1964 Theory of a class of discrete optimal control
systems  {\it J. Electron. Control} {\bf 17}, 697-711

\bibitem[KaMaOr:99]{Kane1}  Kane C Marsden J E and  Ortiz M 1999 Symplectic energy-momentum
integrators {\it J. Math. Phys.} {\bf 40}, 3353-3371


\bibitem[KMOW:00]{Kane3}  Kane C  Marsden J E Ortiz M and  West M 2000 Variational integrators and
the Newmark algorithm for conservative and dissipative mechanical systems
{\it Internat. J. Numer. Math. Eng.} {\bf 49}, 1295-1325.




\bibitem[Lee:83]{Lee1}  Lee T D 1983 Can time be a discrete dynamical variable? {\em Phys. Lett.}, {\bf 122B}, 217-220

\bibitem[Lee:87]{Lee2}  Lee T D 1987 Difference equations and conservation laws {\em J. Statis. Phys.}, {\bf 46}, 843-860


\bibitem[LeMaMD:96]{LMM} de Le\'on M Marrero J C and Mart{\'\i}n de Diego D 1996 Time-dependent constrained Hamiltonian systems and Dirac brackets
{\it Journal of Physics A: Math. Gen.}  {\bf 29} 6843-6859



\bibitem[LeMDSa:02a]{LDA} de Le\'on M Mart{\'\i}n de Diego D and Santamar{\'\i}a A 2002 Geometric integrators and nonholonomic mechanics, {\it Preprint IMAFF-CSIC}

\bibitem[LeMDSa:02b]{LDA1} de Le\'on M, Mart{\'\i}n de Diego D and Santamar{\'\i}a A 2003 Geometric numerical integration of nonholonomic systems and optimal control problems {\it 2nd IFAC Workshop on lagrangian and Hamiltonian Methods for Nonlinear Control, Seville} 2003, 163-168. 



\bibitem[LeoMdD:2002]{LD1} de Le\'on M and Mart{\'\i}n de Diego D 2002 Variational integrators and time-dependent Lagrangian systems {\it Rep. on Math. Phys} {\bf 49} 2/3, 183-192 



\bibitem[Lew:86]{Lewis} Lewis F.L 1986 {\it Optimal Control} (John Wiley\& Sons, New York)


\bibitem[Mae:80]{Maed1} Maeda S 1980 Canonical structure and symmetries for discrete systems {\it Math.
Japonica} {\bf 25}, 405-420

\bibitem[Mae:81]{Maed2} Maeda S 1981 Extension of discrete Noether theorem {\it Math. Japonica} {\bf 26}, 85-90




\bibitem[MaPaSh:98]{Mars2} Marsden J E  Patrick G W  and Shkoller S 1998 Multisymplectic geometry,
variational integrators, and nonlinear PDEs' {\it Comm. Math. Phys.} {\bf 199}, 351-
395


\bibitem[MarWes:01]{Mars6} Marsden J E and  West M 2001 Discrete mechanics and variational integrators {\it  Acta Numerica }, 357-514

\bibitem[MosVes:91]{Mose}  Moser J  and  Veselov A P 1991 Discrete versions of some classical integrable systems and factorization of matrix polynomials {\it Comm. Math. Phys.} {\bf 139}, 
217-243

\bibitem[Nijvan:90]{nv}  Nijmeijer H  and van der Schaft A J 1990  {\it Nonlinear dynamical control systems} (Springer-Verlag, New York)


\bibitem[NorJar:98]{Nort}  Norton K and  Jaroszkiewicz G 1998 Principes of discrete time mechanics, III:
Quantum field theory {\it J. Phys. A} {\bf 31}, 977-1000



\bibitem[PKMO:02]{Pand} Pandolfi A  Kane C  Marsden J E and Ortiz M 2002 Time-discretized variational formulation of nonsmooth frictional contact {\em  Int. J. Num. Methods in
     Engineering} {\bf 53}, 1801-1829


\bibitem[SanCal:94]{Sanz}  Sanz-Serna J M and  Calvo M P 1994 {\it Numerical Hamiltonian Problems} (Chapman\& Hall, London)


\bibitem[Shi:02]{Shi} Shi Y 2002 Symplectic structure of Discrete Hamiltonian Systems,   
{\it J. Math. Anal. Appl.} {\bf 266}, 472-478


\end{thebibliography}
\end{document}